\documentclass{elsart}
\usepackage{amsmath,amssymb}
\usepackage{mathrsfs}
\numberwithin{equation}{section} \allowdisplaybreaks
\newcommand\NN{{\mathbb{N}}}

\def\ls{\leqslant}
\def\gs{\geqslant}
\def\d{{\rm\,d}}
\def\pend{\hfill $\Box$}
\newcommand{\norm}[1]{\left\|#1\right\|}
\newcommand{\abs}[1]{\left|#1\right|}
\newcommand{\set}[1]{\left\{#1\right\}}
\newcommand{\Real}{\mathbb R}
\newcommand{\C}{\mathbb C}
\newcommand{\Z}{\mathbb Z}
\newcommand{\eps}{\varepsilon}
\def\FF{{\mathscr{F}^{-1}}}
\def\F{{\mathscr{F}}}
\def\S{{\mathscr{S}}}
\def\T{{\mathscr{T}}}
\def\M{{\mathscr{M}}}
\def\no{\nonumber}
\def\sgn{\mathrm{sgn}}
\def\supp{\mathrm{supp}}
\begin{document}

\begin{frontmatter}

\title{Wellposedness of Cauchy problem
for the Fourth Order Nonlinear Schr\"odinger Equations in Multi-dimensional Spaces}%

\author[a]{Chengchun Hao\corauthref{*}\thanksref{cchao}},
\ead{hcc@amss.ac.cn}
\author[a]{Ling Hsiao\thanksref{hsiao}},
\ead{hsiaol@amss.ac.cn}
\author[b]{Baoxiang Wang}
\ead{wbx@math.pku.edu.cn} \corauth[*]{Corresponding author.}
 \thanks[cchao]{Supported partially by the Innovation Funds of AMSS, CAS of China.}
\thanks[hsiao]{Supported partially by NSFC(Grant No.10431060), NSAF(Grant No.10276036).}
\address[a]{Academy of Mathematics and Systems
Science, CAS, Beijing 100080, P.R. China}
\address[b]{School of Mathematical
Sciences, Peking University, Beijing 100871, P.R. China}

\date{}
% ----------------------------------------------------------------
\begin{abstract}
 We study the wellposedness of Cauchy problem for the fourth order
nonlinear Schr\"odinger equations
\begin{equation*}
    i\partial_t u=-\eps\Delta u+\Delta^2 u+P\left((\partial_x^\alpha u)_{\abs{\alpha}\ls 2}, (\partial_x^\alpha
    \bar{u})_{\abs{\alpha}\ls 2}\right),\quad t\in \Real,\ x\in\Real^n,
\end{equation*}
where $\eps\in\{-1,0,1\}$, $n\gs 2$ denotes the spatial dimension
and $P(\cdot)$ is a polynomial excluding constant and linear
terms.

\end{abstract}


\end{frontmatter}

% ----------------------------------------------------------------

% ----------------------------------------------------------------
\section{Introduction}

We consider the Cauchy problem for the fourth order nonlinear
Schr\"odinger equations
\begin{align}
    &i\partial_t u=-\eps\Delta u+\Delta^2 u+P((\partial_x^\alpha u)_{\abs{\alpha}\ls 2}, (\partial_x^\alpha
    \bar{u})_{\abs{\alpha}\ls 2}),\quad t\in \Real,\ x\in\Real^n,\label{eq1.1}\\
    &u(0,x)=u_0(x),\label{eq1.2}
\end{align}
where $\eps\in\{-1,0,1\}$, $\alpha\in\Z^n$, $u=u(t,x)$ is a
complex valued wave function of $(t,x)\in\Real^{1+n}$ and $\Delta$
is the Laplace operator on $\Real^n$ with $n\gs 2$. $P(\cdot)$ is
a complex valued polynomial defined in $\C^{n^2+3n+2}$ such that
\begin{align}\label{d2p}
    P(\vec{z})=P(z_1,z_2,\cdots,z_{n^2+3n+2})=\sum\limits_{l\ls\abs{\beta}\ls h\atop \beta\in
    \Z^{n^2+3n+2}}a_\beta z^\beta,\;\text{ for } l,h\in\NN,
\end{align}
and there exists $a_{\beta_0}\neq 0$ for some $\beta_0\in
\Z^{n^2+3n+2}$ with $\abs{\beta_0}=l\gs 2$.

This class of nonlinear Schr\"odinger equation has been widely
applied in many branches in applied science such as deep water wave
dynamics, plasma physics, optical communications and so on
\cite{Dys}. A large amount of interesting work has been devoted to
the study of Cauchy problem to dispersive equations, such as
\cite{Ben,Ca,Chihara,CS,GT,Guo,Hao,Hayashi,Huo,Karpman1996,Kato,KP,KPV1,KPV2,KPV3,KPV98,KPV4,KF84,Pecher,Seg2003,Seg2004,Sj87,Ve88}
and references cited therein.

In order to study the influence of higher order dispersion on
solitary waves, instability and the collapse phenomena, V. I.
Karpman \cite{Karpman1996} introduced a class of nonlinear
Schr\"odinger equations
\begin{align*}
    i\Psi_t+\frac{1}{2}\Delta\Psi+\frac{\gamma}{2}\Delta^2\Psi+f(\abs{\Psi}^2)\Psi=0,\
    x\in\Real^n,\ t\in\Real.
\end{align*}

This system with different nonlinearities was discussed by lots of
authors. In \cite{Seg2003}, by using the method of Fourier
restriction norm, J. Segata studied a special fourth order nonlinear
Schr\"odinger equation in one dimensional space and considered the
three-dimensional motion of an isolated vortex filament, which was
introduced by Da Rios, embedded in inviscid incompressible fluid
fulfilled in an infinite region. And the results have been improved
in \cite{Huo,Seg2004}.

In \cite{Ben}, M. Ben-Artzi, H. Koch and J. C. Saut discussed the
sharp space-time decay properties of fundamental solutions to the
linear equation
\begin{align*}
    i\Psi_t-\varepsilon\Delta\Psi+\Delta^2\Psi=0,\
    \varepsilon\in\{-1,0,1\}.
\end{align*}

In \cite{Guo}, B. L. Guo and B. X. Wang considered the existence and
scattering theory for the Cauchy problem of nonlinear Schr\"odinger
equations with the form
\begin{align}
    &iu_t+(-\Delta)^m u+f(u)=0,\label{w1}\\
    &u(0,x)=\varphi(x),\label{w2}
\end{align}
where $m\gs 1$ is an integer. H. Pecher and W. von Wahl in
\cite{Pecher} proved the existence of classical global solutions
of \eqref{w1}--\eqref{w2} for the space dimensions $n\ls 7m$ for
the case $m\gs 1$. In \cite{Hao}, we discussed the local
wellposedness of the Cauchy problem \eqref{w1}--\eqref{w2} without
smallness of data for $m=2$ and $f(u)=P((\partial_x^j u)_{j\ls 2},
(\partial_x^j \bar{u})_{j\ls 2})$ in one dimension.

In the present paper we deal with the equation \eqref{eq1.1} in
which the difficulty arises from the operator semigroup itself and
the fact that the nonlinearity of $P$ involves the first and the
second derivatives. This could cause the so-called loss of
derivatives so long as we make direct use of the standard methods,
such as the energy estimates, the space-time estimates, etc. In
addition, compared to the one-dimensional case, it is quite
different in the multi-dimensional case, since we have to separate
the spatial space $\Real^n$ into a family of non-overlapping cubes
of bounded size.

  In \cite{KPV3}, C. E. Kenig, G. Ponce and L. Vega  made a great progress on the
nonlinear Schr\"odinger equation of the form
\begin{align*}
    \partial_t u&=i\Delta u+P(u,\bar{u}, \nabla_x u,\nabla_x\bar{u}),\quad t\in\Real,\ x\in\Real^n
\end{align*}
and proved that the corresponding Cauchy problem is locally
wellposed for small data in the Sobolev space $H^s(\Real^n)$ and in
its weighted version by pushing forward the linear estimates
associated with the Schr\"odinger group
$\{e^{it\Delta}\}_{-\infty}^\infty$ and by introducing suitable
function spaces where these estimates act naturally. They also
studied generalized nonlinear Schr\"odinger equations in
\cite{KPV98} and quasi-linear Schr\"odinger equations in
\cite{KPV4}. In the one dimensional case, $n=1$, the smallness
assumption on the size of the data was removed by N. Hayashi and T.
Ozawa \cite{Hayashi} by using a change of variable to obtain an
equivalent system with a nonlinear term independent of $\partial_x
u$, where the new system could be treated by the standard energy
method. And H. Chihara \cite{Chihara}, was able to remove the size
restriction on the data in any dimension by using an invertible
classical pseudo-differential operator of order zero.

  T. Kato showed, in \cite{Kato}, that solutions
of the Korteweg-de Vries equation
\begin{align*}
    \partial_t u+\partial_x^3 u+u\partial_x u=0,\quad t,x\in\Real,
\end{align*}
satisfy that
\begin{align*}
    \int_{-T}^T \int_{-R}^R \abs{\partial_x u(t,x)}^2 \d x\d t\ls
    C(T,R,\norm{u(0,x)}_2).
\end{align*}
And in \cite{KF84}, S. N. Kruzkov and A. V. Faminski{\u\i}, for
$u(0,x) =u_0(x) \in L^2$ such that $x^\alpha u_0(x) \in
L^2((0;+\infty))$, proved that the weak solution of the Korteweg-de
Vries equation constructed there has $l$-continuous space
derivatives for all $t > 0$ if $l < 2\alpha$. The corresponding
version of the above estimate for the Schr\"odinger group
$\{e^{it\Delta}\}_{-\infty}^\infty$
\begin{align*}
    \int_{-T}^T \int_{-R}^R \abs{(1-\Delta)^{1/4} e^{it\Delta}
    u_0}^2 \d x\d t\ls C(T,R,\norm{u_0}_2)
\end{align*}
was simultaneously established by P. Constantin and J.-C. Saut
\cite{CS}, P. Sj\"olin \cite{Sj87} and L. Vega \cite{Ve88} and
others.

In the present paper, we will develop the arguments in \cite{KPV3}
 to study the Cauchy problem \eqref{eq1.1}--\eqref{eq1.2}.
 We first discuss the local smoothing
effects of the unit group $\{S(t)\}_{-\infty}^\infty$ ($S(t)$ is
defined as below) in order to overcome the loss of derivatives in
Sec.2. To construct the work space, we have to study, in Sec.3, the
boundedness properties of the maximal function $\sup_{[0,T]}
\abs{S(t)\cdot}$. This idea is implicit in the splitting argument
introduced by J. Ginibre and Y. Tsutsumi \cite{GT} to deal with
uniqueness for the generalized KdV equation. Finally, we will
consider some special cases in Sec.4 to apply the estimates we have
obtained in the previous sections.

For convenience, we  introduce some notations.
$S(t):=e^{it(\eps\Delta-\Delta^2)}$ denotes the unitary group
generated by $i(\eps\Delta-\Delta^2)$ in $L^2(\Real^n)$. $\bar{z}$
denotes the conjugate of the complex number $z$. $\F u$ or
$\hat{u}$ ($\FF u$, respectively) denotes the Fourier (inverse,
respectively) transform of $u$ with respect to all variables
except the special announcement. $\S$ denotes the space of
Schwartz' functions. $\{Q_\alpha\}_{\alpha\in\Z^n}$ is the family
of non-overlapping cubes of size $R$ such that
$\Real^n=\bigcup\limits_{\alpha\in\Z^n} Q_\alpha$. $D_x^\gamma
f(x)=\FF \abs{\xi}^\gamma\F f(x)$ is the homogeneous derivative of
order $\gamma$ of $f$. $H^s$ ($\dot{H}^s$ , respectively) denotes
the usual inhomogeneous (homogeneous, respectively) Sobolev space.
We denote $\norm{\cdot}_{s,2}=\norm{\cdot}_{H^s}$,
$\norm{\cdot}_p=\norm{\cdot}_{L^p}$ for $1\ls p\ls\infty$ and
\begin{align}\label{nlj}
    \norm{f}_{l, 2, j}=\norm{f}_{H^l(\Real;\, \abs{x}^j\d x)}=\sum_{\abs{\gamma}\ls\,
    l}\left(\int_{\Real^n} \abs{\partial_x^\gamma
    f(x)}^2\abs{x}^j\d x\right)^{1/2}.
\end{align}
 Throughout the paper, the constant $C$ might
be different from each other and $\left[\frac{n}{2}\right]$
denotes the greatest integer that is less than or equal to
$\frac{n}{2}$.

Now we state the main results of this paper.

\begin{thm}[Case $l\gs 3$]\label{thm1}
Let $n\gs 2$. Given any polynomial $P$ as in \eqref{d2p} with $l
\gs 3$, then, for any $u_0\in H^{s}(\Real^n)$ with $s\gs
s_0=n+2+\frac{1}{2}$, there exists a $T=T(\norm{u_0}_{s,2})>0$
such that the Cauchy problem \eqref{eq1.1}-\eqref{eq1.2} has a
unique solution $u(t)$ defined in the time interval $[0,T]$ and
satisfying
\begin{align*}
    u\in C([0,T];H^s(\Real^n)),
\end{align*}
and
\begin{align*}
    u\in L^2([0,T];\, \dot{H}^{s+\frac{1}{2}}(Q)),
\end{align*}
for any cube $Q$ of unit size in $\Real^n$.
\end{thm}

\begin{thm}[Case $l=2$]\label{thm2}
Let $n\gs 2$. Given any polynomial $P$ as in \eqref{d2p} with $l=
2$, then, for any $u_0\in H^{s}(\Real^n)\cap
H^{n+4\left[\frac{n}{2}\right]+8}(\Real^n;\,
\abs{x}^{2\left[\frac{n}{2}\right]+2}\d x)$ with $s\gs
s_0=2n+3\left[\frac{n}{2}\right]+15+\frac{1}{2}$, there exists a
$T=T(\norm{u_0}_{s,2},\,
\norm{u_0}_{n+4\left[\frac{n}{2}\right]+8,2,\left[\frac{n}{2}\right]+2})>0$
such that the Cauchy problem \eqref{eq1.1}-\eqref{eq1.2} has a
unique solution $u(t)$ defined in the time interval $[0,T]$ and
satisfying
\begin{align*}
    u\in C\left([0,T];H^s(\Real^n)\cap
H^{n+\left[\frac{n}{2}\right]+8}(\Real^n;\,
\abs{x}^{2\left[\frac{n}{2}\right]+2}\d x)\right),
\end{align*}
and
\begin{align*}
    u\in L^2([0,T];\, \dot{H}^{s+\frac{1}{2}}(Q)),
\end{align*}
for any cube $Q$ of unit size in $\Real^n$.
\end{thm}

\section{Local Smoothing Effects}

We will prove the local smoothing effects exhibited by the semigroup
$\set{S(t)}_{-\infty}^\infty$ in this section.

\begin{lem}[Local smoothing effect: homogeneous case]
We have the following estimate
\begin{align}
    \sup_{\alpha\in \Z^n}\left(\int_{Q_\alpha}\int_{-\infty}^\infty \abs{D_x^{3/2}S(t) u_0(x)}^2\d t\d x\right)^{1/2}
    \ls C R^{1/2}\norm{u_0}_2, \label{lse1}
    \intertext{and the corresponding dual version}
    \norm{D_x^{3/2}\int_0^T S(-\tau)f(\tau,\cdot)\d\tau}_2
    \ls CR^{1/2}\sum_{\alpha\in \Z^n}\left(\int_{Q_\alpha}\int_0^T\abs{f(x,t)}^2\d t\d x\right)^{1/2}.\label{lse2}
\end{align}
\end{lem}

\textbf{Proof.} \eqref{lse1} and \eqref{lse2} can be derived from
\cite[Theorem 4.1]{KPV1} for which we omit the details.\pend

More precisely, \eqref{lse2} yields, for $t\in[0,T]$, that
\begin{align}
    \norm{D_x^{3/2}\int_0^t S(t-\tau)f(\tau,\cdot)\d\tau}_2
    \ls CR^{1/2}\sum_{\alpha\in \Z^n}\left(\int_{Q_\alpha}\int_0^T\abs{f(x,t)}^2\d t\d x\right)^{1/2}.\label{lse21}
\end{align}

Moreover, with the help of the Fubini theorem, \eqref{lse1} implies
that
\begin{align}
    \sup_{\alpha\in \Z^n}\left(\int_0^T \int_{Q_\alpha}\abs{S(t) u_0(x)}^2\d x\d t\right)^{1/2}
    \ls C R^{1/2}\norm{u_0}_{\dot{H}^{-\frac{3}{2}}}. \label{lse11}
\end{align}
In addition, we have
\begin{align}
    &\sup_{\alpha\in \Z^n}\left(\int_0^T \int_{Q_\alpha}\abs{S(t) u_0(x)}^2\d
    x\d t\right)^{1/2}\no\\
    \ls & T^{\frac{1}{2}}\left(\sup_t \int_{\Real^n}
    \abs{S(t)u_0(x)}^2\d x\right)^{1/2}\no\\
    \ls & T^{\frac{1}{2}}\norm{u_0}_2. \label{lse12}
\end{align}
Interpolating both \eqref{lse11} and \eqref{lse12}, we can obtain
that
\begin{align}\label{lse13}
    \sup_{\alpha\in \Z^n}\left(\int_0^T \int_{Q_\alpha}\abs{S(t) u_0(x)}^2\d
    x\d t\right)^{1/2}\ls
    CR^{\frac{1}{6}}T^{\frac{1}{3}}\norm{u_0}_{\dot{H}^{-\frac{1}{2}}}.
\end{align}

 Now we turn to consider
the inhomogeneous Cauchy problem:
\begin{align}
    i\partial_t u&=-\eps\Delta u+\Delta^2 u+F(t,x),\quad t\in\Real,\, x\in \Real^n,\label{eq2.1}\\
    u(0,x)&=0, \label{eq2.2}
\end{align}
with $F\in\S(\Real\times\Real^n)$. We have the following estimate on
the local smoothing effect in this inhomogeneous case.

\begin{lem}[Local smoothing effect: inhomogeneous
case]\label{lem2.2}
The solution $u(t,x)$ of the Cauchy problem
\eqref{eq2.1}-\eqref{eq2.2} satisfies
\begin{align}\label{lse3}
     \sup_{\alpha\in \Z^n}\norm{D_x^2
    u(t,x)}_{L_x^2(Q_\alpha;\,
    L_t^2([0,T]))}\ls C RT^{\frac{1}{4}}\sum_{\alpha\in \Z^n}\norm{F}_{L_x^2(Q_\alpha;\,
    L_{t}^2([0,T]))}.
\end{align}
\end{lem}

\noindent\textbf{Proof.} Separating
\begin{align*}
F&=\sum_{\alpha\in \Z^n}F \chi_{Q_\alpha}=\sum_{\alpha\in
\Z^n}F_\alpha,
 \intertext{and}
  u&=\sum_{\alpha\in \Z^n}u_\alpha,
\end{align*}
where $u_\alpha(t,x)$ is the corresponding solution of the Cauchy
problem
\begin{align}
    i\partial_t u_\alpha&=-\eps\Delta u_\alpha+\Delta^2 u_\alpha+F_\alpha(t,x),\quad t\in\Real,\, x\in \Real^n,\label{eq21}\\
    u_\alpha(0,x)&=0\label{eq22},
\end{align}
we formally take Fourier transform in both variables $t$ and $x$ in
the equation \eqref{eq21} and obtain
\begin{equation*}
    \hat{u_\alpha}(\tau,\xi)=\frac{\hat{F_\alpha}(\tau,\xi)}{\tau-\eps\abs{\xi}^2-\abs{\xi}^4},\quad
    \textrm{ for each } \alpha\in\Z^n.
\end{equation*}

By the Plancherel theorem in the time variable, we can get
\begin{align}
    &\sup_{\alpha\in \Z^n}\norm{D_x^2
    u_\beta(t,x)}_{L_x^2(Q_\alpha;\
    L_t^2([0,T]))}
    =\sup_{\alpha\in \Z^n}\norm{D_x^2
    \F_t u_\beta(\tau,x)}_{L_x^2(Q_\alpha;\
    L_\tau^2(\Real))}\no\\
    =&\sup_{\alpha\in \Z^n}\norm{D_x^2
    \F_{\xi}^{-1}\left(\frac{\hat{F}_\beta(\tau,\xi)}{\tau-\eps\abs{\xi}^2-\abs{\xi}^4}\right)}_{L_x^2(Q_\alpha;\
    L_\tau^2(\Real))}\no\\
    =& \sup_{\alpha\in \Z^n}\norm{
    \F_{\xi}^{-1}\left(\frac{\abs{\xi}^2}{\tau-\eps\abs{\xi}^2-\abs{\xi}^4}\hat{F}_\beta(\tau,\xi)\right)}_{L_x^2(Q_\alpha;\
    L_\tau^2(\Real))}\no\\
    =& \sup_{\alpha\in \Z^n}\left(\int_{Q_\alpha}\int_{-\infty}^\infty \abs{\F_{\xi}^{-1}
    \left(\frac{\abs{\xi}^2}{\tau-\eps\abs{\xi}^2-\abs{\xi}^4}\hat{F}_\beta(\tau,\xi)\right)}^2\d\tau\d x\right)^{1/2}. \label{eq2.3}
\end{align}

In order to continue the above estimate, we introduce the following
estimate.

\begin{lem}\label{lem2.3}
Let $\M f=\FF m(\xi)\F f$ and
$m(\xi)=\frac{\abs{\xi}^2}{1-\eps\tau^{-\frac{1}{2}}\abs{\xi}^2-\abs{\xi}^4}$
for $\tau>0$ and $\eps\in\{-1,0,1\}$, where $\F(\FF)$ denotes  the
Fourier (inverse, respectively) transform in x only. Then, we have
\begin{align}\label{op}
\sup_{\alpha\in \Z^n} \left(\int_{Q_\alpha}\abs{\M
(g\chi_{Q_\beta})}^2\d x\right)^{1/2}\ls C R
\left(\int_{Q_\beta}\abs{g}^2\d x\right)^{1/2}.
\end{align}
\end{lem}

\noindent\textbf{Proof.} We have
\begin{align*}
    &1-\eps\tau^{-\frac{1}{2}}\abs{\xi}^2-\abs{\xi}^4\\
    &=\left(1-\frac{\eps}{2\sqrt{\tau}}\abs{\xi}^2\right)^2
    -\left(\frac{\eps^2}{4\tau}+1\right)\abs{\xi}^4\\
    &=\left[1+\left(\sqrt{\frac{\eps^2}{4\tau}+1}-\frac{\eps}{2\sqrt{\tau}}\right)\abs{\xi}^2\right]
    \left[1-\left(\sqrt{\frac{\eps^2}{4\tau}+1}+\frac{\eps}{2\sqrt{\tau}}\right)\abs{\xi}^2\right]\\
    &=\left(\frac{2\sqrt{\tau}}{\sqrt{\eps^2+4\tau}-\eps}+\abs{\xi}^2\right)
    \left(\frac{2\sqrt{\tau}}{\sqrt{\eps^2+4\tau}+\eps}-\abs{\xi}^2\right).
\end{align*}

Let $\varphi\in C_0^\infty(\Real)$ with $\supp \varphi\subset
[-1,1]$, $\varphi\equiv 1$ in $[-1/2,1/2]$ and $0\ls \varphi\ls
1$. Denote
$a=\sqrt{\frac{2\sqrt{\tau}}{\sqrt{\eps^2+4\tau}+\eps}}$. Let
$\varphi_1(\xi)=\varphi(2(\abs{\xi}-a))$ and choose
$\varphi_2(\xi)$ such that $\varphi_1(\xi)+\varphi_2(\xi)=1$.

Thus, we have $\supp \varphi_1\subset \{\xi:\,
a-\frac{1}{2}\ls\abs{\xi}\ls a+\frac{1}{2}\}$ and
$\varphi_1(\xi)=1$ for $a-\frac{1}{4}\ls \abs{\xi}\ls
a+\frac{1}{4}$.

Define
\begin{align*}
    \M_j f(x)=\FF m_j(\xi)\F f(x),\ j=1,2,
\end{align*}
where $m_j(\xi)=\varphi_j(\xi)m(\xi)$.

First, we shall establish \eqref{op} for the operator $\M_2$ whose
symbol $m_2(\xi)$ has no singularities. For $p$, $p'$ such that
\begin{align*}
    \frac{1}{p}+\frac{1}{p'}=1 \text{ and }
    \frac{1}{p'}-\frac{1}{p}=\frac{1}{n},
\end{align*}
it follows, by H\"older's inequality, Sobolev's embedding theorem
and the Mihlin theorem,  that
\begin{align*}
    \left(\int_{Q_\alpha}\abs{\M_2(g\chi_{Q_\beta})}^2\d
    x\right)^{\frac{1}{2}}&\ls
    C\norm{1}_{L^{\frac{1}{1/2-1/p}}(Q_\alpha)}
    \norm{\M_2(g\chi_{Q_\beta})}_{L^p(Q_\alpha)}\\
    &\ls C
    R^{n(\frac{1}{2}-\frac{1}{p})}\norm{\M_2(g\chi_{Q_\beta})}_{L^p(Q_\alpha)}\\
    &\ls C
    R^{n(\frac{1}{2}-\frac{1}{p})}\norm{D_x\M_2(g\chi_{Q_\beta})}_{L^{p'}(Q_\alpha)}\\
    &\ls C
    R^{n(\frac{1}{2}-\frac{1}{p})}\norm{g\chi_{Q_\beta}}_{L^{p'}(Q_\alpha)}\\
    &\ls C R\norm{g\chi_{Q_\beta}}_{L^2(Q_\alpha)},
\end{align*}
where we have used the fact
\begin{align}\label{mihlin}
    \norm{D_x \M_2 f}_p\ls C\norm{f}_p, \text{ for } 1<p<\infty.
\end{align}

In fact,
\begin{align*}
    D_x \M_2 f&=\FF \abs{\xi} m_2(\xi)\F f\\
             &=\FF \frac{\abs{\xi}^3[1-
             \varphi(2(\abs{\xi}-a))]}{\left(\frac{2\sqrt{\tau}}{\sqrt{\eps^2+4\tau}-\eps}+\abs{\xi}^2\right)(a^2-\abs{\xi}^2)}\F
             f.
\end{align*}
Denote $\rho(\xi)=\frac{\abs{\xi}^3[1-
             \varphi(2(\abs{\xi}-a))]}{\left(\frac{2\sqrt{\tau}}{\sqrt{\eps^2+4\tau}-\eps}+\abs{\xi}^2\right)(a^2-\abs{\xi}^2)}$.
In the case $\abs{\xi}\gs a+\frac{1}{4}$, we have
\begin{align*}
    \abs{\rho(\xi)}\ls\frac{1-\varphi(2(\abs{\xi}-a))}{4}\ls
    \frac{1}{4}.
\end{align*}
In the case $\abs{\xi}\ls a-\frac{1}{4}$ where we assume that
$a>\frac{1}{4}$, we have the same estimate. Thus,
$\abs{\rho(\xi)}\ls C$. By calculating the derivatives of $\rho$
with respect to the variable $\xi$, we are able to get
\begin{align*}
    \abs{\frac{\partial^\alpha}{\partial \xi^\alpha}\rho(\xi)}\ls
    C\abs{\xi}^{-\abs{\alpha}}, \ \abs{\alpha}\ls L,
\end{align*}
for some integer $L> \frac{n}{2}$. Therefore, we have the desired
estimate \eqref{mihlin} in view of the Mihlin multiplier theorem.

To estimate $\M_1$, we split its symbol $m_1(\xi)$ into a finite
number of pieces (depending only on the dimension $n$). Let
$\theta\in C_0^\infty (\Real)$ with $\supp \theta\subset[-a,a]$
where $a$ is the same as above, define
\begin{align*}
    m_{1,1}(\xi)=m_1(\xi)\theta(4\abs{\bar{\xi}}),\\
    \intertext{and}
    \M_{1,1}f=\FF m_{1,1}(\xi)\F f,
\end{align*}
where $\xi=(\bar{\xi},\xi_n)\in \Real^{n-1}\times\Real$. By a
rotation argument, $m_1(\cdot)$ can be expressed as a finite sum of
$m_{1,1}$'s. Notice that
\begin{align*}
    &\supp m_{1,1}\subset\left\{ \xi=(\bar{\xi},\xi_n)\in
    \Real^{n-1}\times\Real:\ \abs{\bar{\xi}}\ls \frac{a}{4} \text{
    and } a-\frac{1}{2}\ls \abs{\xi}\ls a+\frac{1}{2}\right\},\\
    &Q_\alpha\subset \Real^{n-1}\times [a_\alpha, a_\alpha+R] \text{
    and }    Q_\beta\subset \Real^{n-1}\times [a_\beta, a_\beta+R]
\end{align*}
for appropriate constants $a_\alpha$ and $a_\beta$.

Thus, by the Plancherel theorem in the $\bar{x}$ variable, we have
\begin{align*}
    &\int_{Q_\alpha}\abs{\M_{1,1}(g\chi_{Q_\beta})}^2\!\!\d x\\
    \ls
    &\int_{a_\alpha}^{a_\alpha+R}\!\!\!\int_{\Real^{n-1}}\abs{\M_{1,1}(g\chi_{Q_\beta})}^2\!\!\d\bar{x}\d
    x_n\\
    =
    &\int_{a_\alpha}^{a_\alpha+R}\!\!\!\int_{\Real^{n-1}}
    \abs{\int_{-\infty}^\infty e^{ix_n\xi_n}m_{1,1}(\xi)\widehat{g\chi_{Q_\beta}}(\bar{\xi},\xi_n)\d\xi_n}^2\!\!
    \d\bar{\xi}\d
    x_n\\
    =&\int_{a_\alpha}^{a_\alpha+R}\!\!\!\int_{\Real^{n-1}}
    \abs{\int_{-\infty}^\infty e^{ix_n\xi_n}m_{1,1}(\xi)\int_{a_\beta}^{a_\beta+R}\!\!\!\int_{\Real^{n-1}}
    e^{-iy\xi}(g\chi_{Q_\beta})(\bar{y},y_n)\d\bar{y}\d y_n\d \xi_n}^2\!\!\d\bar{\xi}\d
    x_n\\
    =&\int_{a_\alpha}^{a_\alpha+R}\!\!\!\!\int_{\Real^{n-1}}
    \abs{\int_{a_\beta}^{a_\beta+R}\!\!\!\!\int_{\Real^{n-1}}
    e^{-i\bar{y}\bar{\xi}}(g\chi_{Q_\beta})(\bar{y},y_n)
    \int_{-\infty}^\infty e^{i(x_n-y_n)\xi_n}m_{1,1}(\xi)\d \xi_n\d\bar{y}\d y_n}^2\!\!\d\bar{\xi}\d
    x_n\\
    =&\int_{a_\alpha}^{a_\alpha+R}\!\!\!\!\int_{\Real^{n-1}}
    \abs{\int_{a_\beta}^{a_\beta+R}\!\!\!\!\int_{\Real^{n-1}}
    e^{-i\bar{y}\bar{\xi}}(g\chi_{Q_\beta})(\bar{y},y_n)
    b(x_n,y_n,\bar{\xi})\d\bar{y}\d y_n}^2\!\!\d\bar{\xi}\d
    x_n\\
    =&:E,
\end{align*}
where
\begin{align*}
b(x_n,y_n,\bar{\xi})&=\int_{-\infty}^\infty
e^{i(x_n-y_n)\xi_n}m_{1,1}(\xi)\d \xi_n\\
&=\int_{-\infty}^\infty e^{i(x_n-y_n)\xi_n}\frac{\abs{\xi}^2
\varphi(2(\abs{\xi}-a))}{\left(\frac{2\sqrt{\tau}}{\sqrt{\eps^2+4\tau}-\eps}+\abs{\xi}^2\right)(a^2-\abs{\xi}^2)}
\theta(4\abs{\bar{\xi}})\d \xi_n\\
&=\theta(4\abs{\bar{\xi}})\int_{-\infty}^\infty
e^{i\lambda\xi_n}\frac{\abs{\xi}^2
}{\left(\frac{2\sqrt{\tau}}{\sqrt{\eps^2+4\tau}-\eps}+\abs{\xi}^2\right)(a^2-\abs{\xi}^2)}
\varphi(2(\abs{\xi}-a))\d \xi_n\\
&=K(\lambda,\bar{\xi}),
\end{align*}
here $\lambda=x_n-y_n$ and the support of the integrand is
contained in
\begin{align*}
    A=\left\{ \xi=(\bar{\xi},\xi_n)\in
    \Real^{n-1}\times\Real:\ \abs{\bar{\xi}}\ls \frac{a}{4} \text{
    and } a-\frac{1}{2}\ls \abs{\xi}\ls a+\frac{1}{2}\right\}.
\end{align*}

For $(\bar{\xi},\xi_n)\in A$, we have
$\abs{\xi}^2-a^2=\xi_n^2+\abs{\bar{\xi}}^2-a^2=\xi_n^2-\mu^2$, where
$\mu^2=a^2-\abs{\bar{\xi}}^2>\frac{15}{16}a^2$. Then, we separate
$K$ into two parts, \textit{i.e.}, $K=K_+ +K_-$ with
\begin{align*}
    K_+(\lambda,\bar{\xi})&=\theta(4\abs{\bar{\xi}})\int_0^\infty
e^{i\lambda\xi_n}\frac{1}{\xi_n-\mu}\cdot\frac{\abs{\xi}^2}{\xi_n+\mu}\cdot\frac{
\varphi(2(\abs{\xi}-a))}{\left(\frac{2\sqrt{\tau}}{\sqrt{\eps^2+4\tau}-\eps}+\abs{\xi}^2\right)}
\d \xi_n\\
&=\theta(4\abs{\bar{\xi}})\F_{\xi_n}^{-1}\left(\frac{1}{\xi_n-\mu}\psi(\bar{\xi},\xi_n)\right)(\lambda)\\
&=\theta(4\abs{\bar{\xi}})\F_{\xi_n}^{-1}\left(\frac{1}{\xi_n-\mu}\right)(\lambda)\ast\F_{\xi_n}^{-1}
\psi(\bar{\xi},\xi_n)(\lambda)\\
&=\theta(4\abs{\bar{\xi}})e^{i\lambda\mu}\sgn(\lambda)\ast\F_{\xi_n}^{-1}\left(\psi(\bar{\xi},\xi_n)\right)(\lambda)\\
&=\theta(4\abs{\bar{\xi}})\int_{-\infty}^\infty
e^{iy_n\mu}\sgn(y_n)\F_{\xi_n}^{-1}
\psi(\bar{\xi},\xi_n)(\lambda-y_n)\d y_n,
\end{align*}
where $\psi(\bar{\xi},\xi_n)\in C_0^\infty
(\Real^{n-1}\times\Real^+)$.

Thus, it follows that
\begin{align}\label{K1}
    \abs{K_+(\lambda,\bar{\xi})}\ls\int_{-\infty}^\infty\abs{\F_{\xi_n}^{-1} \psi(\bar{\xi},\xi_n)(y_n)}\d
    y_n\ls C.
\end{align}
In the same way, we have
\begin{align}\label{K2}
    \abs{K_-(\lambda,\bar{\xi})}\ls C.
\end{align}

\eqref{K1} and \eqref{K2} yield that $\abs{b(x_n,
y_n,\bar{\xi})}\ls C$. Now, we can turn to the estimate of
$\int_{Q_\alpha}\abs{\M_{1,1}(g\chi_{Q_\beta})}^2\!\!\d x$ and
get, by the Schwartz inequality in $y_n$ variable and the
Plancherel theorem in $\bar{y}$ variable, that
\begin{align*}
    E&\ls R\int_{a_\alpha}^{a_\alpha+R}\!\!\!\!\int_{\Real^{n-1}}
    \int_{a_\beta}^{a_\beta+R}\abs{\int_{\Real^{n-1}}
    e^{-i\bar{y}\bar{\xi}}(g\chi_{Q_\beta})(\bar{y},y_n)
    \d\bar{y}}^2\!\!b^2(x_n,y_n,\bar{\xi})\d y_n\d\bar{\xi}\d
    x_n\\
    &\ls CR\int_{a_\alpha}^{a_\alpha+R}\!\!\!\!\int_{\Real^{n-1}}
    \int_{a_\beta}^{a_\beta+R}\abs{\int_{\Real^{n-1}}
    e^{-i\bar{y}\bar{\xi}}(g\chi_{Q_\beta})(\bar{y},y_n)
    \d\bar{y}}^2\!\!\d y_n\d\bar{\xi}\d
    x_n\\
    &=CR\int_{a_\alpha}^{a_\alpha+R}\!\!\!\!\int_{a_\beta}^{a_\beta+R}\!\!\!\!\int_{\Real^{n-1}}
    \abs{\int_{\Real^{n-1}}
    e^{-i\bar{y}\bar{\xi}}(g\chi_{Q_\beta})(\bar{y},y_n)
    \d\bar{y}}^2\!\!\d\bar{\xi}\d y_n\d
    x_n\\
    &=CR\int_{a_\alpha}^{a_\alpha+R}\!\!\!\!\int_{a_\beta}^{a_\beta+R}\!\!\!\!\int_{\Real^{n-1}}
    \abs{(g\chi_{Q_\beta})(\bar{y},y_n)
    }^2\!\!\d\bar{y}\d y_n\d x_n\\
    &=CR^2\int_{\Real^n}\abs{(g\chi_{Q_\beta})(y)}^2\!\!\d y.
\end{align*}

Therefore, we obtain
\begin{align*}
    \sup_{\alpha\in\Z^n}\left(\int_{Q_\alpha}\abs{\M(g\chi_{Q_\beta})}^2\d
    x\right)^{1/2}\ls CR\norm{g\chi_{Q_\beta}}_{L^2(\Real^n)}=CR\norm{g}_{L^2(Q_\beta)}
\end{align*}
as desired. \pend

Now, we go back to the proof of Lemma \ref{lem2.2}. We first
consider the part when $\tau>0$ in \eqref{eq2.3}, \textit{i.e.},
\begin{align*}
    &\sup_{\alpha\in \Z^n}\left(\int_{Q_\alpha}\int_0^\infty \abs{\F_{\xi}^{-1}
    \left(\frac{\abs{\xi}^2}{\tau-\eps\abs{\xi}^2-\abs{\xi}^4}\hat{F}_\beta(\tau,\xi)\right)}^2\d\tau\d
    x\right)^{1/2}\\
    =&\left(\int_0^\infty \tau^{\frac{n}{4}-1} \sup_{\alpha\in
    \Z^n}\int_{\tau^{\frac{1}{4}}Q_\alpha}\abs{\int_{\Real^n}
    e^{iy\eta}\frac{\abs{\eta}^2}{1-\eps\tau^{-\frac{1}{2}}\abs{\eta}^2-\abs{\eta}^4}\hat{F}_\beta
    (\tau,\tau^{\frac{1}{4}}\eta)\d\eta}^2 \d y\d\tau\right)^{1/2}\\
    =&\left(\int_0^\infty \tau^{-\frac{n}{4}-1} \sup_{\alpha\in
    \Z^n}\int_{\tau^{\frac{1}{4}}Q_\alpha}\abs{\M \F_t
    F_\beta(\tau,\tau^{-\frac{1}{4}}y)}^2\d y\d\tau\right)^{1/2}\\
    \ls &CR\left(\int_0^\infty \tau^{-\frac{n}{4}-\frac{1}{2}}
    \int_{\tau^{\frac{1}{4}}Q_\beta} \abs{\F_t F_\beta(\tau,\tau^{-\frac{1}{4}}y)}^2\d
    y\d\tau\right)^{1/2}\\
    \ls &CR\left(\int_0^\infty \tau^{-\frac{1}{2}}
    \int_{Q_\beta} \abs{\F_t F(\tau,x)}^2\d
    x\d\tau\right)^{1/2}
\end{align*}
where we have used the changes of variables, Lemma \ref{lem2.3} and
the identity
\begin{align*}
    (\F_\eta^{-1}\hat{F}_\beta(\tau,\tau^{\frac{1}{4}}\eta))(\tau,y)=\tau^{-\frac{n}{4}}\F_t
    F_\beta(\tau,\tau^{-\frac{1}{4}}y).
\end{align*}

For the part when $\tau\in(-\infty, 0)$ in \eqref{eq2.3}, we are
able to split it into two cases. First, it is easier to handle the
case $\eps\in\{0,1\}$ since this  corresponds to the symbol
$\frac{\abs{\eta}^2}{1+\eps\abs{\tau}^{-\frac{1}{2}}\abs{\eta}^2+\abs{\eta}^4}$
which has no singularity. Next, for the case $\eps=-1$, we have
\begin{align*}
    \frac{\abs{\eta}^2}{1+\eps\abs{\tau}^{-\frac{1}{2}}\abs{\eta}^2+\abs{\eta}^4}
    =\frac{\abs{\eta}^2}{\left(1-\frac{1}{2\abs{\tau}^{\frac{1}{2}}}\abs{\eta}^2\right)^2
    +\left(1-\frac{1}{4\abs{\tau}}\right)\abs{\eta}^4}=:B.
\end{align*}
It is obvious that $B$ has no singularity when
$\abs{\tau}>\frac{1}{4}$. Consequently, we have to deal with the
case $\abs{\tau}<\frac{1}{4}$. In fact, we have
\begin{align*}
    B=\frac{\abs{\eta}^2}{\left(\frac{2\abs{\tau}^{\frac{1}{2}}}{1-\sqrt{1-4\abs{\tau}}}-\abs{\eta}^2\right)
      \left(\frac{2\abs{\tau}^{\frac{1}{2}}}{1+\sqrt{1-4\abs{\tau}}}-\abs{\eta}^2\right)
    }.
\end{align*}
Thus, we can use the same argument as in the proof of Lemma
\ref{lem2.3} step by step for the cases $\abs{\eta}^2\ls
\frac{1}{2\sqrt{\abs{\tau}}}$ and $\abs{\eta}^2>
\frac{1}{2\sqrt{\abs{\tau}}}$ respectively, and get that
\begin{align*}
    \sup_{\alpha\in\Z^n}\left(\int_{Q_\alpha}\abs{\F_\eta^{-1} \frac{\abs{\eta}^2}{1-\abs{\tau}^{-\frac{1}{2}}\abs{\eta}^2
    +\abs{\eta}^4}\F_x(g\chi_{Q_\beta})}^2\d
    x\right)^{\frac{1}{2}}\ls CR\left(\int_{Q_\beta}\abs{g}^2\d
    x\right)^{\frac{1}{2}}.
\end{align*}

Therefore,  we obtain, by the Plancherel theorem, the Sobolev
embedding theorem and the H\"older inequality, that
\begin{align*}
   & \sup_{\alpha\in \Z^n}\norm{D_x^2
    u_\beta(t,x)}_{L_x^2(Q_\alpha;\
    L_t^2([0,T]))}\\
    \ls & CR\left(\int_{-\infty}^\infty \abs{\tau}^{-\frac{1}{2}}\int_{Q_\beta} \abs{\F_t F(\tau,x)}^2\d
    x\d\tau\right)^{1/2}\\
    \ls &CR\left(\int_{Q_\beta}\int_{-\infty}^\infty \abs{\abs{\tau}^{-\frac{1}{4}}\F_t F(\tau,x)}^2\d\tau\d
    x\right)^{1/2}\\
    \ls& CR\left(\int_{Q_\beta}\norm{ F(\cdot,x)}_{\dot{H}^{-\frac{1}{4}}([0,T])}^2\d
    x\right)^{1/2}\\
    \ls &CR\left(\int_{Q_\beta}\norm{ F(\cdot,x)}_{L^{\frac{4}{3}}([0,T])}^2\d
    x\right)^{1/2}\\
    \ls& CR\left(\int_{Q_\beta}\norm{1}_{L^{4}([0,T])}^2\norm{ F(\cdot,x)}_{L^2([0,T])}^2\d
    x\right)^{1/2}\\
    \ls & CRT^{\frac{1}{4}}\norm{F}_{L_x^2(Q_\beta;\ L_{t\in[0,T]}^2)},
\end{align*}
which implies the desired result \eqref{lse3}.\pend

\section{Estimates for the Maximal Function}

For simplicity, let $\{Q_\alpha\}_{\alpha\in\Z^n}$ denote the mesh
of dyadic cubes of unit size. We start with  an $L^2$~-continuity
result for the maximal function
$\sup\limits_{[0,T]}\abs{S(t)\cdot}$.

\begin{lem}\label{lem3.1} For any $s>n+\frac{1}{2}$ and $T \in (0,1]$
\begin{align}
    \left(\sum_{\alpha\in\Z^n}\sup_{t\in[0,T]}\sup_{\alpha\in\Z^n}\abs{S(t)u_0(x)}^2\right)^{1/2}\ls
    C\norm{u_0}_{s,2}.\label{m1}\\
    \intertext{In particular,}
    \left(\int_{\Real^n}\sup_{t\in[0,T]}\abs{S(t)u_0(x)}^2\d
    x\right)^{1/2}\ls C\norm{u_0}_{s,2}.\label{m2}
\end{align}
\end{lem}

\noindent\textbf{Proof.}  It is clear that it suffices to prove
the case $t=1$. Let $\{\psi_k\}_{k=0}^\infty$ be a smooth
partition of unity in $\Real^n$ such that the $\psi_k$'s are
radial with $\supp \psi_0\subset\{\xi:\,\abs{\xi}\ls 1\}$, $\supp
\psi_k\subset \{\xi:\, 2^{k-1}\ls\abs{\xi}\ls 2^{k+1}\}$,
$\psi_k(x)\in[0,1]$ and $\abs{\psi_k'(x)}\ls C 2^{-k}$ for
$k\in\Z^+$. For $k\gs 1$ and $t\in[0,2]$, define
\begin{align*}
    I(t,r)=\int_{-\infty}^\infty e^{i\phi_r(s)}\psi_k(s)\d s,
\end{align*}
where the phase function $\phi_r(s)=-\eps ts^2-ts^4+rs$.
Consequently, we have $\phi_r'(s)=-2\eps ts-4ts^3+r$,
$\phi_r''(s)=-2\eps t-12ts^2$ and $\phi_r''''(s)=-24t$.

Denote
\begin{align*}
    &\Omega=\left\{s\in\Real^+:\
    \abs{\phi_r'(s)}<\frac{r}{2}\right\}
    \intertext{and}
    &I_k=[2^{k-1}, 2^{k+1}].
\end{align*}

For $r\in(0,1)$, it is obvious that $\abs{I(t,r)}\ls
\frac{3}{2}\cdot 2^k$.

For $r>1$, we divide it into three cases:

  \ \ i) $\Omega$ is located to the left of $I_k$;

 \  ii)  $\Omega\cap I_k\neq \varnothing$;

   iii)  $\Omega$ is located to the right of $I_k$.

In the cases i) and ii), we have $t\gs C\frac{r}{ 2^{3k}}$ and
$\abs{\phi_r''''(s)}\gs C\frac{r}{ 2^{3k}}$. Thus, with the help of
the Van der Corput lemma, we have $\abs{I(t,r)}\ls
C\left(\frac{2^{3k}}{r}\right)^{1/4}$, where we have used
$\abs{\psi_k'(x)}\ls C 2^{-k}$.

In the case iii), we have $\abs{\phi_r'(s)}\ls \frac{r}{2}$ and
$\frac{r}{t}\ls C2^{3k}$. Integration by parts gives that
\begin{align*}
    \abs{I(t,r)}&\ls \int_0^\infty
    \left[\abs{\frac{\phi_r''(s)}{(\phi_r'(s))^2}}\psi_k(s)+\frac{1}{\abs{\phi_r'(s)}}\abs{\psi_k'(s)}\right]\d
    s\\
    &\ls \int_{I_k}\left[\frac{4(2t+12ts^2)}{r^2}\psi_k(s)
    +\frac{2}{r}\abs{\psi_k'(s)}\right]\d
    s\\
    &\ls C\int_{I_k}\left[\frac{4(2+12\cdot 2^{2(k+1)})}{r\cdot 2^{3k}}\psi_k(s)
    +\frac{2}{r}2^{-k}\right]\d
    s\\
    &\ls \frac{C}{r}\ls C\left(\frac{2^{3k}}{r}\right)^{1/4},
\end{align*}
where we have used the condition $r>1$ in the last step.

Finally, we always have the case in which $\Omega$ is located to
the right of $I_k$ when $r\gs C 2^{3k}$. Hence, we can integrate
by parts as in the case iii) $N$-times and get that
$\abs{I(t,r)}\ls C_N/r^N$ for any $N\in\Z^+$.

Therefore, we have obtained
$$\abs{I(t,r)}\ls \left\{
\begin{array}{ll}
C2^k,&\text{ for } r\in(0,1),\\
C\left(\frac{2^{3k}}{r}\right)^{\frac{1}{4}}, &\text{ for } r\in(1, C2^{3k}),\\
\frac{C_N}{r^N}, &\text{ for } r\gs C2^{3k}.
\end{array}
\right.
$$

In order to continue the proof, we have to introduce some
estimates. The Fourier transform of a radial function
$f(\abs{x})=f(s)$ is given by the formula (\textit{cf.}
\cite[p.154, Theorem 3.3]{Stein})
\begin{align*}
    \hat{f}(r)=\hat{f}(\abs{\xi})=r^{-\frac{n-2}{2}}\int_0^\infty
    f(s)J_{\frac{n-2}{2}}(rs)s^{\frac{n}{2}}\d s,
\end{align*}
and the Bessel function is defined as
\begin{align*}
    J_m(r)=\frac{\left(\frac{r}{2}\right)^m}{\Gamma(m+\frac{1}{2})\Gamma(\frac{1}{2})}
    \int_{-1}^1 e^{irs}(1-s^2)^{\frac{2m-1}{2}}\d s,
\end{align*}
where $\Gamma(\cdot)$ is Gamma functions.

\begin{lem}[\hbox{\cite[Lemma 3.6]{KPV3}}]
\begin{align}
    J_m(r)=O(r^m) &\quad \text{ as } r\to 0,\no\\
    J_m(r)=e^{-ir}\sum_{j=0}^N \alpha_{m,j} r^{-(j+1/2)}
    &+e^{ir}\sum_{j=0}^N \beta_{m,j} r^{-(j+1/2)}\no\\
    &+O(r^{-(N+3/2)}),
    \text{ as } r\to \infty, \label{J1}
\end{align}
for each $N\in\Z^+$.\pend
\end{lem}

We continue the proof of Lemma \ref{lem3.1} now. Denote
\begin{align*}
    \tilde{I}(t,x)&\equiv\abs{\int_{\Real^n}e^{-it(\eps \abs{\xi}^2+\abs{\xi}^4)}e^{ix\xi}\psi_k(\abs{\xi})\d\xi}\\
    &=\frac{1}{r^{(n-2)/2}}\abs{\int_0^\infty
    e^{-it(\eps s^2+s^4)}J_{\frac{n-2}{2}}(rs)\psi_k(s)s^{\frac{n}{2}}\d
    s}\equiv\tilde{I}(t,r),
\end{align*}
where $r=\abs{x}$.

For $r\in(0,1)$, we have
\begin{align*}
    \tilde{I}(t,r)\ls \frac{1}{r^{\frac{n-2}{2}}}\int_0^\infty
    (rs)^{\frac{n-2}{2}}\psi_k(s)s^{\frac{n}{2}}\d s
    \ls C\int_{I_k}s^{n-1}\d s\ls C 2^{nk}.
\end{align*}

For $r>1$, we first consider the remainder term in \eqref{J1} to
obtain the bound
\begin{align*}
    C\frac{1}{r^{\frac{n-2}{2}}}\int_{I_k}
    \frac{1}{(rs)^{N+3/2}}s^{\frac{n}{2}}\d s
    \ls C\frac{2^{kn/2}}{r^n\cdot 2^{kn/2}}\ls \frac{C}{r^n},
\end{align*}
where we have taken N such that $N>\frac{n-1}{2}$.

Next, we deal with the $j$-term in \eqref{J1} for $0\ls j\ls N$
\begin{align*}
    &\frac{1}{r^{\frac{n-2}{2}}}\abs{\int_0^\infty e^{-it(\eps s^2+s^4)}e^{isr}\frac{1}{(sr)^{j+1/2}}s^{n/2}\psi_k(s)\d
    s}\\
    \ls
    &\frac{1}{r^{\frac{n-2}{2}}}\frac{1}{r^{j+1/2}}2^{-(k+1)(j+1/2)}2^{nk/2}
    \abs{\int_0^\infty e^{-it(\eps s^2+s^4)}e^{isr}\psi_k(s)\d s}\\
    \ls&
    C\frac{2^{\frac{kn}{2}-(j+\frac{1}{2})k}}{r^{\frac{n}{2}+j-\frac{1}{2}}}\cdot
    \left\{%
\begin{array}{ll}
    \left(\frac{2^{3k}}{r}\right)^{\frac{1}{4}}, & r\in[1,C2^{3k}], \\
    C_M r^{-M}, &\ r>C2^{3k}.
\end{array}%
\right.
\end{align*}

Now, we consider the $L^1$-norm of $\tilde{I}(t,x)$ with respect to
$x$ variable. For $r\in [0,1]$, it is clear that $\int_{\abs{x}\ls
1}\abs{\tilde{I}(t,x)}\d x\ls C\int_0^1 2^{nk} r^{n-1}\d r\ls C
2^{nk}$. In the case $r\in [1,C2^{3k}]$, we have
\begin{align*}
     \int_{1\ls\abs{x}\ls C2^{3k}}\abs{\tilde{I}(t,x)}\d x
     &\ls C\int_1^{C2^{3k}}
     \frac{2^{\frac{nk}{2}-(j+\frac{1}{2})k}}{r^{\frac{n}{2}+j-\frac{1}{2}}}\left(\frac{2^{3k}}{r}\right)^{\frac{1}{4}}
     r^{n-1}\d r\\
     &\ls
     C2^{\frac{nk}{2}-(j-\frac{1}{4})k}\int_1^{C2^{3k}}r^{\frac{n}{2}-j-\frac{3}{4}}\d
     r
     \ls C2^{k(2n+1)}.
\end{align*}
For $r>C2^{3k}$, it follows that $\int_{\abs{x}\gs
C2^{3k}}\abs{\tilde{I}(t,x)}\d x\ls C2^{kn}$. Thus, we obtain
\begin{align}\label{I}
    \int_{\Real^n}\abs{\tilde{I}(t,x)}\d x\ls C2^{k(2n+1)}.
\end{align}

We decompose the operator $S(t)$ as follows according to the
partition of unity $\{\psi_k\}$,
\begin{align*}
    W_k(t)u_0=\FF e^{-it(\eps \abs{\xi}^2+\abs{\xi}^4)}\psi_k(\xi)\F u_0.
\end{align*}
To prove \eqref{m1}, it suffices to show that
\begin{align}
    &\left(\sum_{\gamma\in\Z^n}\sup_{\abs{t}\ls 1}\sup_{x\in
    Q_\gamma}\abs{\int_{-1}^1
    W_k(t-\tau)g(\tau,\cdot)\d\tau}^2\right)^{1/2}\no\\
    &\ls
    C2^{k(2n+1)}\left(\sum_{\gamma\in\Z^n}\left(\int_{-1}^1\int_{Q_\gamma}\abs{g(t,x)}\d
    x\d t\right)^2\right)^{1/2}.\label{m3}
\end{align}

In fact, we have
\begin{align*}
    &\abs{\int_{-1}^1 W_k(t-\tau)g(\cdot,\tau)\d\tau}\\
    =&\abs{\int_{-1}^1 \FF e^{-i(t-\tau)(\eps \abs{\xi}^2+\abs{\xi}^4)}\psi_k(\xi)\F g(\tau,\xi)\d\tau}\\
    =&\abs{\int_{-1}^1\int_{\Real^n} \left[\FF e^{-i(t-\tau)(\eps \abs{\xi}^2+\abs{\xi}^4)}\psi_k(\xi)\right](y)
    \cdot g(\tau,x-y)\d y\d\tau}\\
    =&\abs{\int_{\Real^n}\int_{-1}^1 \left[\FF e^{-i(t-\tau)(\eps \abs{\xi}^2+\abs{\xi}^4)}\psi_k(\xi)\right](y)
    \cdot g(\tau,x-y)\d\tau\d y}\\
    \ls&\int_{\Real^n}\int_{-1}^1 \left[\abs{\FF
    e^{-i(t-\tau)(\eps \abs{\xi}^2+\abs{\xi}^4)}\psi_k(\xi)\right](y)}
    \abs{ g(\tau,x-y)}\d\tau\d y\\
    \ls&\int_{\Real^n}\sup_{t\in [0,2]}\tilde{I}(t,\abs{y})\int_{-1}^1\abs{g(\tau,x-y)}\d\tau\d y\\
    \ls&\sum_{\alpha\in\Z^n}\left(\sup_{t\in [0,2]}\sup_{y\in
    Q_\alpha}\tilde{I}(t,\abs{y})\right)\int_{-1}^1\int_{Q_\alpha}\abs{g(\tau,x-y)}\d
    y\d\tau.
\end{align*}
Thus, the left hand side of \eqref{m3} is bounded by
\begin{align}\label{m4}
    \left(\sum_{\gamma\in\Z^n}\left[\sum_{\alpha\in\Z^n}\left(\sup_{t\in [0,2]}\sup_{y\in
    Q_\alpha}\tilde{I}(t,\abs{y})\right)\sup_{x\in Q_\gamma}\int_{-1}^1\int_{Q_\alpha}\abs{g(\tau,x-y)}\d
    y\d\tau\right]^2\right)^{1/2}.
\end{align}
Let $E_{\alpha,\gamma}:=Q_\alpha-Q_\gamma=2^nQ_\alpha-x_\gamma$,
where $x_\gamma$ denotes the center of $Q_\gamma$. Then, we can get
\begin{align*}
    \sup_{x\in Q_\gamma}\int_{-1}^1\int_{Q_\alpha}\abs{g(\tau,x-y)}\d
    y\d\tau\ls\int_{-1}^1\int_{E_{\alpha,\gamma}}\abs{g(\tau,z)}\d
    z\d\tau,
\end{align*}
and consequently \eqref{m4} is bounded, with the help of the
Minkowski inequality, by
\begin{align*}
    &\left(\sum_{\gamma\in\Z^n}\left[\sum_{\alpha\in\Z^n}\left(\sup_{t\in
    [0,2]}\sup_{y\in
    Q_\alpha}\tilde{I}(t,\abs{y})\right)\int_{-1}^1\int_{E_{\alpha,\gamma}}\abs{g(\tau,z)}\d
    z\d\tau  \right]^2\right)^{1/2}\\
    \ls &\sum_{\alpha\in\Z^n}\left[\left(\sup_{t\in
    [0,2]}\sup_{y\in
    Q_\alpha}\tilde{I}(t,\abs{y})\right)\left(\sum_{\gamma\in\Z^n}\left(\int_{-1}^1\int_{E_{\alpha,\gamma}}\abs{g(\tau,z)}\d
    z\d\tau  \right)^2\right)^{1/2}\right]\\
    \ls &C\sum_{\alpha\in\Z^n}\left[\left(\sup_{t\in
    [0,2]}\sup_{y\in
    Q_\alpha}\tilde{I}(t,\abs{y})\right)\left(\sum_{\gamma\in\Z^n}\left(\int_{-1}^1\int_{Q_\gamma}\abs{g(\tau,z)}\d
    z\d\tau  \right)^2\right)^{1/2}\right],
\end{align*}
where we have used the fact that $E_{\alpha,\gamma}$ can be
covered by a finite number (independent of $\alpha$) of $Q_\gamma$
in the last step.

From \eqref{I}, we have the estimate
\begin{align*}
    \sum_{\alpha\in\Z^n}\sup_{t\in
    [0,2]}\sup_{y\in
    Q_\alpha}\tilde{I}(t,\abs{y})=\sup_{t\in
    [0,2]}\int_{\Real^n}\tilde{I}(t,\abs{y})\d y\ls C2^{(2n+1)k}.
\end{align*}

 Therefore, the proof is completed. \pend

When we deal with the Cauchy problem \eqref{eq1.1}-\eqref{eq1.2} in
the case $l=2$, we shall use the $l^1$-estimate of the maximal
function $\sup_{[0,T]}\abs{S(t)u_0(x)}$ as the following inequality.

\begin{lem}\label{lem3}
We have the estimate
\begin{align}
    &\sum_{\alpha\in\Z^n}\sup_{[0,T]}\sup_{x\in
    Q_\alpha}\abs{S(t)u_0(x)}\no\\
    \ls & C(1+T)^{\left[\frac{n}{2}\right]+2}\left(\norm{u_0}_{n+3\left[\frac{n}{2}\right]+7,\, 2}
    +\norm{u_0}_{n+3\left[\frac{n}{2}\right]+7,\, 2,\, 2\left[\frac{n}{2}\right]+2}\right),\label{m5}
\end{align}
where $\norm{\cdot}_{l,\, 2,\, j}$ is defined as in \eqref{nlj}.
\end{lem}

\noindent\textbf{Proof.}  Taking $t_0\in [0,T]$ such that
\begin{align*}
    \abs{f(t_0)}\ls \frac{\int_0^T\abs{f(t)}\d t}{T},
\end{align*}
we have for any $t$
\begin{align*}
    f(t)=f(t_0)+f(t)-f(t_0)=f(t_0)+\int_{t_0}^t f'(s)\d s.
\end{align*}
Thus, we can get
\begin{align*}
    \abs{f(t)}&\ls \frac{1}{T}\int_0^T\abs{f(t)}\d
    t+\int_{t_0}^t\abs{f'(s)}\d s\\
    &\ls \frac{1}{T}\int_0^T\abs{f(t)}\d
    t+\int_0^T\abs{\partial_t f(t)}\d t,
\end{align*}
namely,
\begin{align}\label{eq3}
    \sup_{t\in[0,T]} \abs{f(t)}\ls \frac{1}{T}\int_0^T\abs{f(t)}\d
    t+\int_0^T\abs{\partial_t f(t)}\d t.
\end{align}

Notice that
\begin{align}\label{eq4}
    \norm{f}_{L^1(\Real^n)}\ls
    C\norm{f}_{L^2(\Real^n)}+C\norm{\abs{x}^{\bar{n}}f}_{L^2(\Real^n)},
\end{align}
where $\bar{n}>n/2$, and
\begin{align}
    &x_jS(t)f=S(t)(x_j f)+4itS(t)(\partial_{x_j}\Delta f)-2\eps it S(t)(\partial_{x_j}f),\label{eq5}\\
    x_kx_jS(t)f=&((-2\eps it)^2+8it)S(t)(\partial_{x_k}\partial_{x_j}f)+4itS(t)(\delta_{kj}\Delta
    f)-2\eps itS(t)(\delta_{kj}f)\no\\
    &\quad +16\eps t^2 S(t)(\partial_{x_k}\partial_{x_j}\Delta
    f)+(4it)^2S(t)(\partial_{x_k}\partial_{x_j}\Delta^2
    f)\no\\
    &\quad -2\eps itS(t)(x_k\partial_{x_j}f+x_j\partial_{x_k}f)+4itS(t)(x_k\partial_{x_j}\Delta f+x_j\partial_{x_k}\Delta
    f)\no\\
    &\quad +S(t)(x_k x_j f),\label{eq6}
\end{align}
where $\delta_{kj}=1$ if $k=j$, $\delta_{kj}=0$ if $k\neq j$.

We have, from \eqref{eq3}-\eqref{eq6}, the Sobolev embedding theorem
and the Fubini theorem, that
\begin{align*}
    &\sum_{\alpha\in\Z^n}\sup_{[0,T]}\sup_{x\in
    Q_\alpha}\abs{S(t)u_0(x)}\\
    \ls &\sum_{\alpha\in\Z^n}\sup_{[0,T]}\left(\int_{Q_\alpha}\abs{S(t)u_0(x)}\d
    x+\sum_{\abs{\beta}\ls n}\int_{Q_\alpha}\abs{\partial_x^\beta S(t)u_0(x)}\d
    x\right)\\
    \ls&\sum_{\alpha\in\Z^n}\left(\int_{Q_\alpha}\sup_{[0,T]}\abs{S(t)u_0(x)}\d
    x+\sum_{\abs{\beta}\ls n}\int_{Q_\alpha}\sup_{[0,T]}\abs{\partial_x^\beta S(t)u_0(x)}\d
    x\right)\\
    \ls
    &\sum_{\alpha\in\Z^n}\Big\{\frac{C}{T}\int_{Q_\alpha}\int_0^T
    \abs{S(t)u_0(x)}\d t\d x+\int_{Q_\alpha}\int_0^T
    \abs{\partial_t S(t)u_0(x)}\d t\d x\\
    &\qquad +\sum_{\abs{\beta}\ls n}\int_{Q_\alpha}\big(\frac{C}{T}\int_0^T
    \abs{S(t)\partial_x^\beta u_0(x)}\d t+\int_0^T
    \abs{\partial_t S(t)\partial_x^\beta u_0(x)}\d t\big)\d
    x\Big\}\\
    \ls&\frac{C}{T}\int_{\Real^n}\int_0^T\abs{S(t)u_0(x)}\d t\d x+
    \int_{\Real^n}\int_0^T\abs{S(t)\Delta u_0(x)}\d t\d x\\
    & +\int_{\Real^n}\int_0^T\abs{S(t)\Delta^2 u_0(x)}\d t\d
    x +\sum_{\abs{\beta}\ls n}\frac{C}{T}\int_{\Real^n}\int_0^T\abs{S(t)\partial_x^\beta u_0(x)}\d
    t\d x\\
    &+\sum_{\abs{\beta}\ls n}\int_{\Real^n}\int_0^T\abs{S(t)\partial_x^\beta \Delta u_0(x)}\d
    t\d x+\sum_{\abs{\beta}\ls n}\int_{\Real^n}\int_0^T\abs{S(t)\partial_x^\beta \Delta^2 u_0(x)}\d
    t\d x\\
    \ls&\frac{C}{T}\int_0^T\int_{\Real^n}\abs{S(t)u_0(x)}\d x\d t
    +\frac{C}{T}\sum_{\abs{\beta}\ls n}\int_0^T\int_{\Real^n}\abs{S(t)\partial_x^\beta u_0(x)}\d x\d
    t\\
    &\qquad\qquad +C\sum_{\abs{\beta}\ls n+4}\int_0^T\int_{\Real^n}\abs{S(t)\partial_x^\beta u_0(x)}\d x\d
    t\\
    \ls&\frac{C}{T}\int_0^T[\norm{S(t)u_0}_2+\norm{\abs{x}^{\bar{n}}S(t)u_0}_2]\d
    t\\
    &\qquad\qquad+\frac{C}{T}\sum_{\abs{\beta}\ls n}\int_0^T[\norm{S(t)\partial_x^\beta u_0}_2
    +\norm{\abs{x}^{\bar{n}}S(t)\partial_x^\beta u_0}_2]\d
    t\\
    &\qquad\qquad\qquad+C\sum_{\abs{\beta}\ls n+4}\int_0^T[\norm{S(t)\partial_x^\beta u_0}_2
    +\norm{\abs{x}^{\bar{n}}S(t)\partial_x^\beta u_0}_2]\d
    t\\
    \ls&
    C(1+T)^{\bar{n}+1}(\norm{u_0}_{n+3\bar{n}+4,\,2}+\norm{u_0}_{n+3\bar{n}+4,\,2,\,2\bar{n}}),
\end{align*}
which implies the desired result if we choose
$\bar{n}=\left[\frac{n}{2}\right]+1$. \pend

\section{The Wellposedness}

We will give the proofs of the main Theorems in this section.

\noindent\textbf{Proof of Theorem \ref{thm1}.} Similar to the proof
in \cite[Theorem 4.1]{KPV3}, we shall only consider the most
interesting case $s=s_0$. The general case follows by combining this
result with the fact that the highest derivatives involved in that
proof always appear linearly and with some commutator estimates
(\textit{see} \cite{KP}) for the cases where $s\neq k+1/2$,
$k\in\Z^+$. For simplicity of the exposition, we shall assume that
\begin{align*}
    P((\partial_x^\alpha u)_{\abs{\alpha}\ls 2}, (\partial_x^\alpha
    \bar{u})_{\abs{\alpha}\ls 2})=\partial_{x_j}^2 u\partial_{x_k}^2
    u\partial_{x_m}^2 u
\end{align*}

For $u_0\in H^{s_0}(\Real^n)$, we denote by $u=\T(v)=\T_{u_0}(v)$
the solution of the linear inhomogeneous Cauchy problem
\begin{align}
    i\partial_t u&=-\eps\Delta u+\Delta^2 u+\partial_{x_j}^2 v\partial_{x_k}^2
    v\partial_{x_m}^2 v,\quad t\in \Real, x\in \Real^n,\label{d31}\\
    u(0,x)&=u_0(x),\label{d32}
\end{align}
In order to construct $\T$ being a contraction mapping in some
space, we use the integral equation
\begin{align}\label{inteq}
    u(t)=\T (v)(t)=S(t)u_0-i\int_0^t S(t-\tau)\partial_{x_j}^2 v(\tau)\partial_{x_k}^2
    v(\tau)\partial_{x_m}^2 v(\tau)\d\tau.
\end{align}

We introduce the following work space
\begin{align*}
    Z_T^E=\Big\{w:[0,T]\times&\Real^n\to\C:\ \sup_{t\in[0,T]}\norm{w(t)}_{s_0,2}\ls  E; \\
    &\sum_{\abs{\beta}=s_0+1/2}\sup_{\alpha\in \Z^n}\left(\int_0^T \int_{Q_\alpha}\abs{\partial_x^\beta w(t,x)}^2
    \d x\d t\right)^{1/2}\ls T^\delta;\\
    &\left(\sum_{\alpha\in \Z^n}\sup_{t\in[0,T]}\sup_{x\in Q_\alpha}\abs{D_x^2 w(t,x)}^2
    \right)^{1/2}\ls E\Big\},
\end{align*}
where $\delta<1/3$ is a constant.

We notice, for any $\beta\in\Z^n$ with $\abs{\beta}=s_0-3/2$, that
\begin{align*}
    \partial_x^\beta (\partial_{x_j}^2 v\partial_{x_k}^2
    v\partial_{x_m}^2 v)=&\partial_x^\beta \partial_{x_j}^2 v\partial_{x_k}^2
    v\partial_{x_m}^2 v+ \partial_{x_j}^2 v\partial_x^\beta\partial_{x_k}^2
    v\partial_{x_m}^2 v+ \partial_{x_j}^2 v\partial_{x_k}^2
    v\partial_x^\beta\partial_{x_m}^2 v\\
    &+R((\partial_x^\gamma v)_{2\ls \abs{\gamma}\ls
    s_0-1/2}),
\end{align*}
where $\abs{\beta'}=\abs{\beta}-1$ and $\abs{\beta''}=1$ with
$\beta'$, $\beta''\in\Z^n$.

 Thus, from the integral equation \eqref{inteq},
 \eqref{lse13}, \eqref{lse1}, the Sobolev embedding theorem and the commutator estimates \cite{KP}, we can get, for
$v\in Z_T^D$, that
\begin{align}
    &\sum_{\abs{\beta}=s_0+1/2}\sup_{\alpha\in \Z^n}\left(\int_0^T \int_{Q_\alpha}\abs{\partial_x^\beta u(t,x)}^2
    \d x\d t\right)^{1/2}\no\\
    \ls &\sum_{\abs{\beta}=s_0+1/2}\sup_{\alpha\in \Z^n}\left(\int_0^T \int_{Q_\alpha}\abs{S(t)\partial_x^\beta u_0(x)}^2
    \d x\d t\right)^{1/2}\no\\
    &\qquad+\sum_{\abs{\beta}=s_0+1/2}\sup_{\alpha\in \Z^n}
    \left(\int_0^T \int_{Q_\alpha}\abs{\int_0^t S(t-\tau)\partial_x^\beta (\partial_{x_j}^2 v\partial_{x_k}^2
    v\partial_{x_m}^2 v)}^2\d x\d t\right)^{1/2}\no\\
    \ls & CT^{\frac{1}{3}}\norm{u_0}_{s_0,
    2}+CT^{\frac{1}{4}}\sum_{\abs{\beta}=s_0-3/2}\sum_{j,k,m=1}^n
    \sum_{\alpha\in\Z^n} \norm{\partial_x^\beta \partial_{x_j}^2 v\partial_{x_k}^2
    v\partial_{x_m}^2 v}_{L_x^2(Q_\alpha;\, L_t^2([0,T]))}\no\\
    &\qquad +C\int_0^T\norm{D_x^{\frac{1}{2}}R((\partial_x^\gamma v)_{2\ls \abs{\gamma}\ls
    s_0-1/2})}_2\d t\no\\
    \ls & CT^{\frac{1}{3}}\norm{u_0}_{s_0, 2}+
    CT^{\frac{1}{4}}\sum_{\abs{\beta}=s_0+1/2}\sup_{\alpha\in \Z^n}\left(\int_0^T \int_{Q_\alpha}\abs{\partial_x^\beta v(t,x)}^2
    \d x\d t\right)^{1/2}\no\\
    & \qquad\cdot\left(\sum_{\alpha\in\Z^n}\sup_{t\in[0,T]}\sup_{x\in Q_\alpha}\abs{D_x^2
    v(t,x)}^2\right)+CT\sup_{t\in[0,T]}\norm{v}_{s_0,\, 2}^3\no\\
    \ls & CT^{\frac{1}{3}}\norm{u_0}_{s_0, 2}+CT^{\frac{1}{4}}T^\delta E^2+CT
    E^3\no\\
    \ls &T^\delta,\label{c1}
\end{align}
where we have taken $T$ so small that
\begin{align}\label{Tc1}
    CT^{\frac{1}{3}-\delta}\norm{u_0}_{s_0, 2}+ CT^{\frac{1}{4}}E^2+CT^{1-\delta} E^3\ls 1,
\end{align}
in the last step.

By the Sobolev embedding theorem, \eqref{lse2} and the H\"older
inequality, we have
\begin{align}
    \sup_{t\in[0,T]}\norm{u(t)}_{s_0,2}\ls
    &\norm{u_0}_{s_0,2}+\sup_{t\in[0,T]}\int_0^t \norm{S(t-\tau)\partial_{x_j}^2 v(\tau)\partial_{x_k}^2
    v(\tau)\partial_{x_m}^2 v(\tau)}_2\d\tau\no\\
    & \quad +\sup_{t\in[0,T]}\norm{D_x^{\frac{3}{2}}\int_0^t S(t-\tau)
    D_x^{s_0-\frac{3}{2}} \partial_{x_j}^2 v(\tau)\partial_{x_k}^2
    v(\tau)\partial_{x_m}^2 v(\tau)\d\tau}_2\no\\
    \ls &\norm{u_0}_{s_0,2}+T\sup_{t\in[0,T]}\norm{\partial_{x_j}^2 v(t)\partial_{x_k}^2
    v(t)\partial_{x_m}^2 v(t)}_2\no\\
    &\quad+\sum_{\alpha\in\Z^n}\left(\int_{Q_\alpha}\int_0^T\abs{D_x^{s_0-\frac{3}{2}}
    (\partial_{x_j}^2 v(t)\partial_{x_k}^2
    v(t)\partial_{x_m}^2 v(t))}^2\d
    t\d x\right)^{1/2}\no\\
    \ls &\norm{u_0}_{s_0,2}+T\sup_{t\in[0,T]}\norm{v}_{\frac{n}{3}+2,2}^3\no\\
    &\quad +\sum_{j,k,m=1}^n \sum_{\alpha\in\Z^n}\left(\int_{Q_\alpha}\int_0^T\abs{D_x^{s_0-\frac{3}{2}}
    \partial_{x_j}^2 v(t)\partial_{x_k}^2 v(t)\partial_{x_m}^2 v(t)}^2\d
    t\d x\right)^{1/2}\no\\
    &\quad+\sum_{\alpha\in\Z^n}\left(\int_{Q_\alpha}\int_0^T\abs{(D_x^\gamma v)_{2\ls\abs{\gamma}\ls s_0-\frac{1}{2}}}^2\d
    t\d x\right)^{1/2}\no\\
    \ls &\norm{u_0}_{s_0,2}+T\sup_{t\in[0,T]}\norm{v}_{\frac{n}{3}+2,2}^3\no\\
    &\quad
    +\sum_{\abs{\beta}=s_0+\frac{1}{2}}\sup_{\alpha\in\Z^n}\left(\int_0^T\int_{Q_\alpha}\abs{\partial_x^\beta
    v}^2\d x\d t\right)^{1/2}\cdot \sum_{\alpha\in\Z^n}\sup_{t\in[0,T]}\sup_{x\in
    Q_\alpha}\abs{D_x^2 v}^2\no\\
    &\quad+T^{\frac{1}{2}}\sup_{t\in[0,T]}\norm{v}_{s_0,2}^3\no\\
    \ls&\norm{u_0}_{s_0,2}+(T+T^{\frac{1}{2}})E^3+T^\delta E^2\no\\
    \ls& E,\label{c2}
\end{align}
where in the last step, we have chosen
 $T$ small enough such that
\begin{align}\label{Tc2}
    (T+T^{\frac{1}{2}})E^2+T^\delta E\ls \frac{1}{2}.
\end{align}

Similar to the derivation of \eqref{c1} and \eqref{c2}, we obtain,
by inserting \eqref{m2} in \eqref{inteq}, that
\begin{align}
    &\left(\sum_{\alpha\in \Z^n}\sup_{t\in[0,T]}\sup_{x\in Q_\alpha}\abs{D_x^2 u(t,x)}^2
    \right)^{1/2}\no\\
    \ls& C\norm{u_0}_{s_0,2}+C\left(\sum_{\alpha\in \Z^n}\sup_{t\in[0,T]}\sup_{x\in Q_\alpha}
    \abs{\int_0^t S(t-\tau) D_x^2 \partial_{x_j}^2 v(\tau)\partial_{x_k}^2 v(\tau)\partial_{x_m}^2 v(\tau)}^2
    \right)^{1/2}\no\\
    \ls &C\norm{u_0}_{s_0,2}+
    CT\sup_{t\in[0,T]}\norm{v}_{s_0,2}^3\no\\
    \ls &\norm{u_0}_{s_0,2}+CTE^3\no\\
    \ls &E,\label{c3}
\end{align}
where we have taken
\begin{align}\label{Ec}
    E=2C\norm{u_0}_{s_0,2},
\end{align}
and $T$ sufficiently small such that
\begin{align}\label{Tc3}
    CTE^2\ls 1.
\end{align}

Therefore, choosing an $E$ as in \eqref{Ec} and then taking a $T$
sufficiently small such that \eqref{Tc1}, \eqref{Tc2} and
\eqref{Tc3} hold, we obtain that the mapping
\begin{align*}
    \T=\T_{u_0}: Z_T^E\to Z_T^E
\end{align*}
is well defined.

For convenience, we denote
\begin{align*}
    \Lambda_T(w)=\max&\left\{ET^{-\delta}\sum_{\abs{\beta}=s_0+1/2}\sup_{\alpha\in \Z^n}
    \left(\int_0^T \int_{Q_\alpha}\abs{\partial_x^\beta w(t,x)}^2
    \d x\d t\right)^{1/2} ;\right.\no\\
     &\left.\sup_{t\in[0,T]}\norm{w(t)}_{s_0,2};\,
    \left(\sum_{\alpha\in \Z^n}\sup_{t\in[0,T]}\sup_{x\in Q_\alpha}\abs{D_x^2 w(t,x)}^2
    \right)^{1/2} \right\}.
\end{align*}

To show that $\T$ is a contraction mapping, we apply the estimates
obtained in \eqref{c1}, \eqref{c2} and \eqref{c3} to the following
integral equation
\begin{align*}
    \T v(t)-\T w(t)=\int_0^t S(t-\tau) &\big[\partial_{x_j}^2 v(\tau)\partial_{x_k}^2 v(\tau)\partial_{x_m}^2 v(\tau)\\
    &-\partial_{x_j}^2 w(\tau)\partial_{x_k}^2 w(\tau)\partial_{x_m}^2 w(\tau)\big](\tau)\d\tau,
\end{align*}
and obtain, for $v,w \in Z_T^E$, that
\begin{align*}
    \Lambda_T(\T v-\T w)&\ls
    CT^\delta\Lambda_T(v-w)\cdot[\Lambda_T^2(v)+\Lambda_T^2(w)]\no\\
     &\ls 2CT^\delta E^2\Lambda_T(v-w).
\end{align*}
where the constant $C$ depends only on the form of $P(\cdot)$ and
the linear estimates \eqref{lse1}, \eqref{lse21}, \eqref{lse3} and
\eqref{m1}.

Thus, we can choose $0<T\ll 1$ satisfying \eqref{Tc1}, \eqref{Tc2},
\eqref{Tc3} and
\begin{align}\label{Tc4}
    2CT^\delta E^2\ls 1/2.
\end{align}

Therefore, for those $T$, satisfying \eqref{Tc1}, \eqref{Tc2},
\eqref{Tc3} and \eqref{Tc4}, the mapping $\T_{u_0}$ is a contraction
mapping in $Z_T^E$. Consequently, by the Banach contraction mapping
principle, there exists a unique function $u\in Z_T^E$ such that
$\T_{u_0}u=u$ which solves the Cauchy problem.

By the method given in \cite[Theorem 4.1]{KPV3}, we can prove the
persistence property of $u(t)$ in $H^{s_0}$, \textit{i.e.}
\begin{align*}
    u(t,x)\in C\left([0,T];H^{s_0}(\Real^n)\right),
\end{align*}
the uniqueness and the continuous dependence on the initial data of
solution. For simplicity, we omit the rest of the proof.\pend

\noindent\textbf{Proof of Theorem \ref{thm2}.} For simplicity, we
assume
\begin{align*}
    P((\partial_x^\alpha u)_{\abs{\alpha}\ls 2}, (\partial_x^\alpha
    \bar{u})_{\abs{\alpha}\ls 2})=\abs{\Delta u}^2.
\end{align*}

It will be clear, from the argument presented below, that this
does not represent any loss of generality. And as in the proof of
Theorem \ref{thm1}, we consider the case
$s=s_0=2n+3\left[\frac{n}{2}\right]+15+\frac{1}{2}$.

We introduce the following work space
\begin{align*}
    X_T^E=\Big\{ w: [0,T]\times\Real^n \to \C:\
    &\sup_t\norm{w(t)}_{s_0,2}\ls E;\
    \sup_t\norm{w(t)}_{n+\left[\frac{n}{2}\right]+8, 2, 2\left[\frac{n}{2}\right]+2}\ls E;\\
    &\sum_{\abs{\beta}=s_0+1/2}\sup_{\alpha\in \Z^n}\left(\int_0^T \int_{Q_\alpha}\abs{\partial_x^\beta w(t,x)}^2
    \d x\d t\right)^{1/2}\ls T^\delta;\\
    (1+T)^{-\left[\frac{n}{2}\right]-2}&\left(\sum_{\alpha\in \Z^n}\sup_{t\in[0,T]}\sup_{x\in Q_\alpha}\abs{D_x^2 w(t,x)}^2
    \right)^{1/2}\ls E\Big\},
\end{align*}
where $0<\delta<1/4$ is a constant and the norm is defined as
\begin{align*}
    \norm{v(t)}_{X_T^E}=\max\Big\{&\sup_t\norm{v(t)}_{\frac{23}{2}};\
    \sup_t\norm{w(t)}_{n+\left[\frac{n}{2}\right]+8, 2,
    2\left[\frac{n}{2}\right]+2};\\
    &ET^{-\delta}\sum_{\abs{\beta}=s_0+1/2}\sup_{\alpha\in \Z^n}\left(\int_0^T \int_{Q_\alpha}
    \abs{\partial_x^\beta w(t,x)}^2 \d x\d t\right)^{1/2};\\
    &(1+T)^{-\left[\frac{n}{2}\right]-2}
    \left(\sum_{\alpha\in \Z^n}\sup_{t\in[0,T]}\sup_{x\in Q_\alpha}\abs{D_x^2 w(t,x)}^2
    \right)^{1/2}
    \Big\}.
\end{align*}

For $u_0\in H^{s_0}(\Real^n)\cap
H^{n+4\left[\frac{n}{2}\right]+8}(\Real^n;\,
\abs{x}^{2\left[\frac{n}{2}\right]+2}\d x)$ and $v\in X_T^E$, we
denote by $u=\T(v)=\T_{u_0}(v)$ the solution of the linear
inhomogeneous Cauchy problem
\begin{align}
    i\partial_t u&=-\eps\Delta u+\Delta^2 u+\abs{\Delta v}^2, \label{l21}\\
    u(0,x)&=u_0(x),\label{l22}
\end{align}
and consider the corresponding integral equation
\begin{align}\label{l23}
    u(t)=\T(v)(t)=S(t)u_0-i\int_0^t S(t-\tau)\abs{\Delta
    v}^2\d\tau.
\end{align}

We notice that
\begin{align*}
    \partial_x^\beta (\Delta v \Delta\bar{v})=\partial_x^\beta \Delta v
    \Delta\bar{v}+\Delta
    v\partial_x^\beta \Delta\bar{v}+R_0\left((\partial_x^{\gamma_1}v
    \partial_x^{\gamma_2}v)_{\abs{\gamma_1},\, \abs{\gamma_2}\ls
    s_0-\frac{1}{2}}\right),
\end{align*}
where $\beta$, $\gamma_1$, $\gamma_2\in\Z^n$.

From the integral equation \eqref{l23}, \eqref{lse13}, the Sobolev
embedding theorem and \eqref{lse3},  we can get, as in \eqref{c1},
that
\begin{align}
    &\sum_{\abs{\beta}=s_0+\frac{1}{2}}\sup_{\alpha\in\Z^n}\left(\int_0^T\int_{Q_\alpha}\abs{\partial_x^\beta
    u(t,x)}^2 \d x\d t\right)^{1/2}\no\\
    \ls &\sum_{\abs{\beta}=s_0+\frac{1}{2}}\sup_{\alpha\in\Z^n}\left(\int_0^T\int_{Q_\alpha}\abs{\partial_x^\beta
    S(t)u_0(x)}^2 \d x\d t\right)^{1/2}\no\\
    &\qquad +\sum_{\abs{\beta}=s_0+\frac{1}{2}}\sup_{\alpha\in\Z^n}\left(\int_0^T\int_{Q_\alpha}\abs{
    \int_0^tS(t-\tau)\partial_x^\beta\abs{\Delta v(\tau)}^2\d\tau}^2 \d x\d
    t\right)^{1/2}\no\\
    \ls &C
    T^{\frac{1}{3}}\norm{u_0}_{s_0,2}+CT^{\frac{1}{4}}\!\!\!\!\sum_{\abs{\beta}=s_0-\frac{3}{2}}\sum_{\alpha\in\Z^n}
    \left(\int_0^T\int_{Q_\alpha}\abs{\partial_x^\beta \Delta v
    \Delta\bar{v}}^2\d x\d t\right)^{1/2}\no\\
    &\qquad +C\int_0^T\norm{D_x^{\frac{1}{2}}R_0}_2\d t\no\\
    \ls &CT^{\frac{1}{3}}\norm{u_0}_{s_0,2}+CT^{\frac{1}{4}}\!\!\!\!\sum_{\abs{\beta}=s_0+\frac{1}{2}}
    \sum_{\alpha\in\Z^n} \left(\int_0^T\int_{Q_\alpha}\abs{\partial_x^\beta v}^2\d x\d
    t\right)^{1/2}\no\\
    &\qquad \cdot \left(\sum_{\alpha\in\Z^n}\sup_t\sup_{x\in
    Q_\alpha}\abs{D_x^2 v}\right)+CT\sup_t \norm{v}_{s_0,2}^2\no\\
    \ls &C
    T^{\frac{1}{3}}\norm{u_0}_{s_0,2}+CT^{\frac{1}{4}}T^\delta
    E+CTE^2\no\\
    \ls T^\delta,\label{lc1}
\end{align}
if we take $T$ sufficiently small such that
\begin{align}\label{lT1}
    CT^{\frac{1}{3}-\delta}\norm{u_0}_{s_0,2}+CT^{\frac{1}{4}}E+CT^{1-\delta}E^2\ls
    1.
\end{align}

We can rewrite \eqref{l23} as
\begin{align}
    u(t)=S(t)\left(u_0-i\int_0^t S(-\tau)\abs{\Delta v(\tau)}^2\d\tau\right).
\end{align}

From \eqref{m5}, we have
\begin{align*}
    &(1+T)^{-\left[\frac{n}{2}\right]-2}
    \sum_{\alpha\in\Z^n}\sup_t\sup_{x\in Q_\alpha}\abs{D_x^2
    u(t,x)}\\
    \ls&  C(\norm{u_0}_{n+3\left[\frac{n}{2}\right]+9,2}
    +\norm{u_0}_{n+3\left[\frac{n}{2}\right]+9,2,2\left[\frac{n}{2}\right]+2})\\
     &\qquad +\sum_{\alpha\in\Z^n}\sup_t\sup_{x\in Q_\alpha}\abs{D_x^2
    S(t)\int_0^t S(-\tau)\abs{\Delta v(\tau)}^2\d\tau},
\end{align*}
where the second term in the right hand side can be bounded by
\begin{align*}
   & C\sup_t \norm{\int_0^t S(-\tau)\abs{\Delta v}^2\d\tau}_{n+3\left[\frac{n}{2}\right]+9,2}
    +\sup_t\norm{\int_0^t S(-\tau)\abs{\Delta v}^2\d\tau}_{n+3\left[\frac{n}{2}\right]+9,2,
    2\left[\frac{n}{2}\right]+2}\\
    \ls &C \int_0^T \norm{\abs{\Delta
    v}^2}_{n+3\left[\frac{n}{2}\right]+9,2}\d\tau +C\int_0^T \norm{S(-\tau)\abs{\Delta
    v}^2}_{n+3\left[\frac{n}{2}\right]+9,2,
    2\left[\frac{n}{2}\right]+2}\d\tau\\
    \ls & CT\sup_t
    \norm{v}_{n+3\left[\frac{n}{2}\right]+11,2}^2+CT\sup_t \sum_{\abs{\gamma}\ls n+3\left[\frac{n}{2}\right]+9}
    \norm{\abs{x}^{\left[\frac{n}{2}\right]+1}S(-\tau)\partial_x^\gamma \abs{\Delta
    v}^2}_2\\
    \ls & CT\sup_t \norm{v}_{n+3\left[\frac{n}{2}\right]+11,2}^2
    +CT\left\{\sup_t \sum_{\abs{\gamma}\ls n+3\left[\frac{n}{2}\right]+9}
    \norm{\abs{x}^{\left[\frac{n}{2}\right]+1}\partial_x^\gamma \abs{\Delta
    v}^2}_2\right.\\
    &\quad\left.+\sum_{\abs{\gamma}\ls n+6\left[\frac{n}{2}\right]+12}t^{\left[\frac{n}{2}\right]+1}
    \norm{\partial_x^\gamma \abs{\Delta v}^2}_2
    +C\sum_{\abs{\beta}+\abs{\gamma}\ls n+6\left[\frac{n}{2}\right]+10 \atop 1\ls\abs{\beta}\ls\left[\frac{n}{2}\right]}
    t^{\abs{\beta}} \norm{x^\beta \partial_x^\gamma \abs{\Delta
    v}^2}_2\right\}\\
    \ls &CT\sup_t \norm{v}_{n+3\left[\frac{n}{2}\right]+11,2}^2
    + CT\sup_t \sum_{\abs{\gamma}\ls n+3\left[\frac{n}{2}\right]+9}
    \sum_{\abs{\gamma_1}+\abs{\gamma_2}=\abs{\gamma}+4 \atop \abs{\gamma_1}\ls \abs{\gamma_2}}
    \norm{\partial_x^{\gamma_2}v}_\infty
    \norm{\abs{x}^{\left[\frac{n}{2}\right]+1}\partial_x^{\gamma_1}v}_2\\
    &\quad +CT^{\left[\frac{n}{2}\right]+2}\sup_t
    \norm{v}_{n+6\left[\frac{n}{2}\right]+14,2}\\
    &\quad+\sum_{\abs{\beta}+\abs{\gamma}\ls
    n+4\left[\frac{n}{2}\right]+12
    \atop 1\ls \abs{\beta}\ls  \left[\frac{n}{2}\right]}T^{\abs{\beta}+1}\sup_t
    \sum_{\abs{\gamma_1}+\abs{\gamma_2}=\abs{\gamma}+4 \atop \abs{\gamma_1}\ls \abs{\gamma_2}}
    \norm{\partial_x^{\gamma_2}v}_\infty
    \norm{x^\beta\partial_x^{\gamma_1}v}_2\\
    \ls& CT\sup_t \norm{v}_{n+3\left[\frac{n}{2}\right]+11,2}^2
    +CT\sup_t
    \norm{v}_{n+4\left[\frac{n}{2}\right]+14,2}\sup_t\norm{v}_{n+\left[\frac{n}{4}\right]+7,2,2\left[\frac{n}{2}\right]+2}\\
    &\quad +CT^{\left[\frac{n}{2}\right]+2}\sup_t\norm{v}_{n+6\left[\frac{n}{2}\right]+14,2}\\
    &\quad
    +C\sum_{1\ls\abs{\beta}\ls\left[\frac{n}{2}\right]}T^{\abs{\beta}+1}\sup_t
    \norm{v}_{\frac{3n}{2}+4\left[\frac{n}{2}\right]+15+\frac{1}{2},2}\sup_t
    \norm{v}_{n+\left[\frac{n}{2}\right]+8,2,2\left[\frac{n}{2}\right]+2}\\
    \ls &CT\sup_t\norm{v}_{s_0,2}^2+CT\sup_t\norm{v}_{s_0,2}\sup_t\norm{v}_{n+\left[\frac{n}{2}\right]+8,2,2\left[\frac{n}{2}\right]+2}
    +CT^{\left[\frac{n}{2}\right]+2}\sup_t\norm{v}_{s_0,2}^2\\
    &\quad
    +CT(T+T^{\left[\frac{n}{2}\right]})\sup_t\norm{v}_{s_0,2}\sup_t\norm{v}_{n+\left[\frac{n}{2}\right]+8,2,2\left[\frac{n}{2}\right]+2}\\
    \ls & CT(1+T^{\left[\frac{n}{2}\right]+1})E^2\\
    \ls &\frac{1}{2}E,
\end{align*}
if we take $E\gs 2C(\norm{u_0}_{n+3\left[\frac{n}{2}\right]+9,2}
    +\norm{u_0}_{n+4\left[\frac{n}{2}\right]+8,2,2\left[\frac{n}{2}\right]+2})$ and $T$ so small that
\begin{align}\label{lT2}
    CT(1+T^{\left[\frac{n}{2}\right]+1})E\ls \frac{1}{2}.
\end{align}

Thus, we obtain
\begin{align}\label{lc2}
    (1+T)^{-\left[\frac{n}{2}\right]-2}
    \sum_{\alpha\in\Z^n}\sup_t\sup_{x\in Q_\alpha}\abs{D_x^2 u(t,x)} \ls E.
\end{align}

Similar to \eqref{lc2}, we have
\begin{align}
   & \sup_t \norm{u(t)}_{n+\left[\frac{n}{2}\right]+8,2,2\left[\frac{n}{2}\right]+2}
    =\sup_t \sum_{\abs{\gamma}\ls
    n+\left[\frac{n}{2}\right]+8}\norm{\abs{x}^{\left[\frac{n}{2}\right]+1}\partial_x^\gamma
    u(t)}_2\no\\
    \ls &\sum_{\abs{\gamma}\ls n+\left[\frac{n}{2}\right]+8}
    \left[\sup_t \norm{\abs{x}^{\left[\frac{n}{2}\right]+1}S(t)\partial_x^\gamma
    u_0}_2+\sup_t \norm{\int_0^t\abs{x}^{\left[\frac{n}{2}\right]+1}S(t-\tau)\partial_x^\gamma \abs{\Delta v}^2\d\tau}_2
    \right]\no\\
    \ls &\sup_t\sum_{\abs{\gamma}\ls n+\left[\frac{n}{2}\right]+8}
    \norm{\abs{x}^{\left[\frac{n}{2}\right]+1}\partial_x^\gamma u_0}_2
    +C\sup_t \sum_{\abs{\gamma}\ls
    n+4\left[\frac{n}{2}\right]+11}t^{\left[\frac{n}{2}\right]+1}\norm{\partial_x^\gamma
    u_0}_2\no\\
    &\quad+ C\sup_t\sum_{\abs{\beta}+\abs{\gamma}\ls n+4\left[\frac{n}{2}\right]+9 \atop 1\ls \abs{\beta}\ls\left[\frac{n}{2}\right]}
    t^{\abs{\beta}}\norm{x^\beta \partial_x^\gamma u_0}_2
    +\sup_t \int_0^t \sum_{\abs{\gamma}\ls n+\left[\frac{n}{2}\right]+8}
    \norm{\abs{x}^{\left[\frac{n}{2}\right]+1}\partial_x^\gamma \abs{\Delta
    v}^2}_2\no\\
    &\quad +C\sup_t \sum_{\abs{\gamma}\ls
    n+4\left[\frac{n}{2}\right]+11}t^{\left[\frac{n}{2}\right]+1}\norm{\partial_x^\gamma
    \abs{\Delta v}^2}_2+C\sup_t\sum_{\abs{\beta}+\abs{\gamma}\ls
    n+4\left[\frac{n}{2}\right]+9
    \atop 1\ls \abs{\beta}\ls\left[\frac{n}{2}\right]}
    t^{\abs{\beta}}\norm{x^\beta \partial_x^\gamma \abs{\Delta
    v}^2}_2\no\\
    \ls &
    C(1+T)^{\left[\frac{n}{2}\right]+1}(\norm{u_0}_{n+4\left[\frac{n}{2}\right]+8,2,2\left[\frac{n}{2}\right]+2}
    +\norm{u_0}_{s_0,2})\no\\
    &\quad +T\sup_t\norm{v}_{s_0,2}\norm{v}_{n+\left[\frac{n}{2}\right]+8,2,
    2\left[\frac{n}{2}\right]+2}+
    CT^{\left[\frac{n}{2}\right]+1}\norm{v}_{s_0,2}^2\no\\
    \ls &E. \label{lc3}
\end{align}

As in the derivation of \eqref{c2}, we can get
\begin{align}
    \sup_t\norm{u(t)}_{s_0,2}\ls E.\label{lc4}
\end{align}

The estimates \eqref{lc1}, \eqref{lc2}, \eqref{lc3} and \eqref{lc4}
yield $u\in X_T^E$. Thus, we can fix $E$ and $T$ as above such that
$\T$ is a contraction mapping from $X_T^E$ to itself. By the
standard argument used in the proof of Theorem \ref{thm1}, we can
complete the proof for which we omit the details. \pend

\section*{Acknowledgment}

The authors would like to thank the referees for valuable comments
and suggestions for the original manuscript.

% ----------------------------------------------------------------
%\bibliographystyle{amsplain}
%\bibliography{}

\end{document}